\documentclass[12pt]{article}

\usepackage[english]{babel}
\usepackage{amsfonts}
\usepackage{amssymb}
\usepackage{amsmath}
\usepackage{amsthm}
\usepackage{enumerate}

\theoremstyle{plain}
\newtheorem*{mainthm}{Main Theorem}

\theoremstyle{plain}
\newtheorem{thm}{Theorem}[section]
\newtheorem{prop}[thm]{Proposition}
\newtheorem{lem}[thm]{Lemma}
\newtheorem{cor}[thm]{Corollary}
\newtheorem{defn}[thm]{Definition}

\newtheorem{example}[thm]{Example}

\theoremstyle{remark}
\newtheorem*{rem}{Remark}

\newcommand{\set}[1]{\left\{#1\right\}}

\newcommand{\bR}{{\mathbb{R}}}
\newcommand{\bN}{{\mathbb{N}}}
\newcommand{\bZ}{{\mathbb{Z}}}

\title{A family of meta-{F}ibonacci sequences\\ defined by variable order
recursions}

\author{Nathaniel D. Emerson}

\begin{document}

\maketitle

\begin{abstract}
We define a family of meta-Fibonacci sequences where the order of
the of recursion at stage $n$ is a variable $r(n)$, and the
$n^{th}$ term of a sequence is the sum of the previous $r(n)$
terms.  For the terms of any such sequence, we give upper and
lower bounds which depend only on $r(n)$.
\end{abstract}

\begin{center}
California State University Channel Islands\\
One University Drive\\
Camarillo, California 93012-8599\\
Email: \texttt{Nathaniel.Emerson@csuci.edu}

\end{center}

\section{Introduction}

We consider \emph{meta-Fibonacci sequences}, that is a sequence
given by a Fibonacci-type recursion, where the recursion varies
with the index.  We describe a new family of meta-Fibonacci
sequences defined by variable order recursions and give
closed-from upper and lower bounds for the terms of any such
sequence. This family is considerably different from previously
described families of meta-Fibonacci sequences (see \cite{DGNW},
\cite{HT93b}, \cite{CCT} and \cite{JR05}) both in terms of its
definition and behavior. See \cite{CCT} for a nice history of the
subject.

In this paper we denote integer valued sequences by Roman letters,
and other sequences by Greek letters.

 We define a family of meta-Fibonacci sequences.  The regular Fibonacci numbers are of
course obtained by adding the previous two terms of a sequence:
$f_n = f_{n-1} + f_{n-2}$.  If we add the previous three terms, we
obtain the \emph{Tribonacci numbers}: $t_n = t_{n-1} + t_{n-2} +
t_{n-3}$. If we add the previous $r$ terms we obtain the
\emph{$r$-generalized Fibonacci numbers} (the ``$r$-bonacci
numbers''): $f_{r,n} = f_{r,n-1} + \dots + f_{r,n-r}$. Now let $r$
vary as a function of $n$.  We call the resulting numbers
\emph{variable-$r$ meta-Fibonacci numbers.}  Let $\bN$ denote the
non-negative integers and $\bZ^+$ denote the positive integers.

\begin{defn} \label{defn: vrmFs}

Let $r: \bN \to \bZ^+$ such that $r(0) = 1$ and $r(n)$ is
\emph{sublinear}, that is $ r(n) \leq n$ for all $n \geq 1$.
Define
\[
    b(n) = \sum_{k=1}^{r(n)} b(n-k), \quad n > 1,
\]
with initial condition $b(0) =1$.  We call the sequence $b(n)$ a
\emph{variable-$r$ meta-Fibonacci sequence,} and say that $r(n)$
\emph{generates} $b(n)$.
\end{defn}

For brevity, we call $b(n)$ a \emph{variable-$r$-bonacci
sequence}. It is clear that any such sequence is a non-decreasing
sequence of positive integers. Additionally, it is clear distinct
$r(n)$ generate distinct sequences $b(n)$. So we have defined an
uncountable family of meta-Fibonacci sequences, in one-to-one
correspondence with sublinear sequences of positive integers.  In
this paper, we examine the dependence of $b(n)$ on $r(n)$.

\begin{example}
If $r(1)= 1$ and $r(n) = 2$ for all $n \geq 2$, then $b(n)  =
f_{n+1}$, where $(f_n)$ is the usual Fibonacci sequence.
\end{example}

\begin{example}
  If $r(n) =1$ for all $n$, then
$b(n)=1 $ for all $n$.  If $r(n) =1$ for all $n$ large, then
$b(n)$ is eventually constant.
\end{example}

\begin{example} \label{ex: b(n) max growth}
 If $r(n) = n$ for $n \geq 1$, then $b(n) = 2^{n-1}$ for $n \geq 1$.

\[
\begin{array}{|c|c|c|c|c|c|c|c|c|c|c|}\hline
  n     & 0 & 1 &  2 &  3 & 4  & 5  &  6 & 7  & 8  & 9  \\
  \hline
  r(n) &  1 &  1 &  2 &  3 &  4 & 5  &  6 &  7 & 8  & 9  \\
  \hline
  b(n) & 1  &  1 & 2  & 4  & 8  &  16 & 32  & 64  & 128  & 256  \\
  \hline
\end{array}
\]

\end{example}

We give an estimate for $b(n)$, which depends only on $n$, $r(1)$,
\dots, $r(n)$.  That is, a \emph{closed-form} estimate. We define
two quantities which we will use to estimate $b(n)$.

\begin{defn}
Let $b(n)$ be a variable-$r$ meta-Fibonacci sequence generated by
$r(n)$. For $n \geq 1$ define
\[
    \lambda(n) = 1 + \frac{r(n) - 1}{r(n-1)}.
\]
\end{defn}

\begin{defn}
Let $b(n)$ be a variable-$r$ meta-Fibonacci sequence generated by
$r(n)$. For $n \geq s \geq 1$ define
\[
    \mu(n,s) = 2 + [ r(n) - r(n-1) - 1]\prod_{k=n-s}^{n-1} 1/r(k) .
\]
\end{defn}

We use $\lambda(n)$ and $\mu(n,s)$ to estimate the growth rate of
$b(n)$, that is the ratio of successive terms. We obtain
closed-from upper and lower bounds.

\begin{mainthm} \label{thm: growth est}
Let $b(n)$ be a variable-$r$ meta-Fibonacci sequence generated by
$r(n)$. For all $n \geq 1$
\[
    \min \set{\lambda(n), \mu(n,r(n)- 1)} \leq  \frac{b(n)}{b(n-1)} \leq \max \set{\lambda(n),
    \mu(n, r(n-1)) }.
\]
\end{mainthm}

We will prove the Main Theorem in $\S$\ref{sect: Growth}.   We
obtain the following bounds on $b(n)$.

\begin{cor} \label{cor: b(n) est}
Let $b(n)$ be a variable-$r$ meta-Fibonacci sequence generated by
$r(n)$. For all $n \geq 2$ we have
\[
   \prod_{k=2}^{n} \min \set{\lambda(k), \mu(k,r(k)- 1)} \leq  {b(n)}
   \leq \max  \prod_{k=2}^{n} \set{\lambda(k), \mu(k, r(k-1)) }.
\]
\begin{proof}
Note that $r(1) = 1$, so $b(1) = 1$ and $b(n) = b(n)/b(1)$.  Write
$b(n)$ as a telescoping product,
\[
    b(n) = \prod_{k=2}^{n} \frac{b(k)}{b({k-1})},
\]
and apply the Main Theorem.
\end{proof}
\end{cor}

Variable-$r$ meta-Fibonacci sequences are considerably different
than any meta-Fibonacci sequence that the author is familiar with.
 Hofstadter's $Q$-sequence (\cite[p. 137]{Hof79} and \cite{Cnly89}) may be taken as a typical
example.  Let $Q(1) = Q(2) = 1$ and
\[
    Q(n) = Q(n-Q(n-1)) + Q(n-Q(n-2)), \quad n > 2.
\]
The recursion for $Q(n)$ is ``self-referential'' and the order of
the recursion is fixed. The terms we add to obtain $Q(n)$ are not
necessarily the immediately previous terms. So for some $n$, we
may add terms that are early in the sequence (and thus small) and
$Q(n)$ will be small. In contrast for a sequence $b(n)$, an
``external'' variable $r(n)$ controls the recursion, and the order
of the recursion is not generally fixed.  We always add the
immediately previous $r(n)$ terms, so $b(n)$ is always the sum of
the largest of the previous terms. These differences result in
$Q(n)$ having much more complicated behavior than $b(n)$.  The
behavior of $Q(n)$ has been described as ``chaotic,'' \cite{Hof79}
while $b(n)$ is non-decreasing.  They also result in very
different rates of growth for the two type of sequences, see
$\S$\ref{sect: Asymptotics}.

 Variable-$r$-bonacci numbers
were originally discovered by the author while studying dynamical
systems \cite{E01}; specifically the dynamics of complex
polynomials.  In the study of dynamical systems, one can consider
\emph{closest return times}--most intuitively, the iterates of a
given point under some map that are closer to the point than any
previous iterate. In \cite{E05} it is shown that certain
generalized closest return times of polynomials are
\emph{extended} variable-$r$ meta-Fibonacci numbers (see
Definition \ref{defn: ex vrmFs}).  This result generalizes the
fact that there exist polynomials whose closest return times are
the ordinary Fibonacci sequence \cite[Ex. 12.4]{BH92}.

The remainder of this paper is organized as follows. In
$\S$\ref{sect: Growth} we prove the Main Theorem.  We give a
series of estimates, and then combine them.   We show that the
growth rate of any variable-$r$-bonacci sequence is at most
exponential.  For variable-$r$-bonacci sequences with $r(n) > 1$
for all $n$ sufficiently large, we show that the growth rate is at
least exponential. In $\S$\ref{sect: Asymptotics} we the study the
asymptotics of $b(n)$. We show that wide variety of growth rates
occur--exponential, linear, and logarithmic. In contrast, for many
other meta-Fibonacci sequences the growth rate is linear. In
$\S$\ref{sect: Generalization} We define a generalization of
$b(n)$, which is defined for all integers and removes restrictions
on $r$.

\section{Estimates on Growth} \label{sect: Growth}

Consider the growth rate of the ordinary Fibonacci numbers, that
is the ratio of successive Fibonacci numbers.  It is well known
that the growth rate is exponential, and converges to the Golden
Section: $(1 + \sqrt{5})/2$.  By a similar argument, the growth
rate of the $r$-bonacci numbers converges to $ \alpha_r$, the
unique real root of the polynomial
\[
    x^r - x^{r-1} - \dots -x -1
\]
with $1 < \alpha_r < 2$, (all other roots have complex modulus
less than 1), see \cite{Miles}.  In this section, we examine the
growth rate of variable $r$-bonacci sequences.  We give a series
of estimates on the growth of variable-$r$-bonacci numbers, which
we will combine to prove our Main Theorem. Throughout this
section, let $b(n)$ be a variable-$r$ meta-Fibonacci sequence
generated by $r(n)$

We derive basic information about the limiting behavior of $b(n)$.

\begin{lem}
The sequence $b(n)$ is eventually constant if and only if \newline
$\limsup_{n \to \infty} r(n) = 1$.
\end{lem}

\begin{lem}
We have $\lim_{n \to \infty} b(n) = \infty$ if and only if
$\limsup_{n \to \infty} r(n) > 1$.
\end{lem}

Thus, a variable-$r$-bonacci sequence converges if and only if it
is eventually constant. . Clearly for a given $n$, the larger
$r(n)$ is, the larger $b(n)$ will be. However in many of these
estimates, it is $\Delta r(n) = r(n) - r(n-1)$ which most strongly
influences the growth rate. The following lemma is the our basic
estimate; we give a condition for $b(n)$ to double.

\begin{lem} \label{lem: Delta r(n) = 1 -> growth =2}
If $\Delta r(n) = 1$ for some $n \geq 1$, then $b(n)/b(n-1) =2$.
\begin{proof}
We have $r(n) = r(n-1) + 1$ for some $n$. Hence
\begin{align*}
     b(n) &= \sum_{k=1}^{r(n)} b(n-k)\\
        &= b(n-1) + \sum_{k=2}^{r(n-1)+1} b(n-k)\\
        &=  b(n-1) + \sum_{i=1}^{r(n-1)} b({n-1 - i})\\
        &= 2b(n-1).\\
\end{align*}
\end{proof}
\end{lem}

We extend the above lemma to cover all cases for $\Delta r(n)$.
 We obtain fairly complete information on the short-term growth of
 $b(n)$, particularly on the relative magnitude of $b(n)/b(n-1)$ and 2.

\begin{thm}\label{thm: b(n) growth}
For all $n \geq 1$ the following hold:

\begin{enumerate}[\indent a.]
    \item $b(n) / b(n-1) = 1 $ if and only if $\Delta r(n) =
    1-r(n-1)$.
    \item $1 < b(n) / b(n-1) < 2 $ if and only if $1-r(n-1) < \Delta r(n) <  1$.
    \item $b(n) / b(n-1) = 2 $ if and only if $\Delta r(n) =  1$.
    \item $ b(n) / b(n-1) > 2 $ if and only if $\Delta r(n) >  1$.
\end{enumerate}

\begin{proof}
We will prove the ``if'' part of each case.  Case \textit{a} is
equivalent to $r(n) =1$, so it is clear. Case \textit{c} is Lemma
\ref{lem: Delta r(n) = 1 -> growth =2}.  In case \textit{b} we
have $r(n) < r(n-1) + 1$, so
\[
    b(n) = \sum_{k=1}^{r(n)} b(n-k) < \sum_{k=1}^{r(n-1)+1}
    b(n-k)= 2 b(n-1),
\]
by Lemma \ref{lem: Delta r(n) = 1 -> growth =2}.  Case \textit{d}
is similar.  The ``only if'' directions follow by considering the
above cases.

\end{proof}
\end{thm}

We can start to examine the long-term behavior of the growth of
$b(n)$.

\begin{cor} \label{cor: lim sup growth < 2 -> r(n) const}
If $\displaystyle \limsup_{n \to \infty} \frac{b(n)}{b(n-1)} < 2$,
then $r(n)$ is eventually constant.
\begin{proof}
By Theorem \ref{thm: b(n) growth}.b for all $n$ sufficiently large
$r(n)  - 1 < r(n-1) $.  It follows that for $n$ large, $r(n)$ is a
non-increasing sequence of positive integers, so is eventually
constant.
\end{proof}
\end{cor}

We give a universal upper bound for $b(n)$, one which does not
depend on $r(n)$.

\begin{lem} \label{lem: max growth 2 n-1}
For all $n \geq 1 $, we have $b(n) \leq 2^{n-1}$.
\begin{proof}
 Let $\hat{r}(n) = n$ for all $n \geq 1$ and let
$\hat{b}(n)$ be the variable-$r$-bonacci sequence generated by
$\hat{r}(n)$.  By Lemma \ref{lem: Delta r(n) = 1 -> growth =2},
$\hat{b}(n) = 2^{n-1}$ for all $n \geq 1$.  Note that $b(0) =
\hat{b}(0) = 1$ and inductively
\[
    b(n) = \sum_{k=1}^{r(n)} b(n-k) \leq \sum_{k=1}^{n} b(n-k) \leq
    \sum_{k=1}^{\hat{r}(n)} \hat{b}({n-k}) = \hat{b}(n) = 2^{n-1}.
\]
\end{proof}
\end{lem}

This bound shows that all variable-$r$-bonacci sequences are
$O(2^{n-1})$.  That is, at worst exponential order.

The following lemma is the basis for many of our other estimates.
We relate the growth of $b(n)$ to $r(n)$.

\begin{lem} \label{lem: b(n) <= r(n)}
For all $n \geq 1 $,
\[
    \frac{b(n)}{b(n-1)} \leq r(n).
\]
\begin{proof}
\begin{align*}
    b(n) &= \sum_{k=1}^{r(n)} b(n-k) \\
        &\leq \sum_{k=1}^{r(n)} b(n-1) \qquad \text{since the $b(n)$ are
        non-increasing,}\\
        &= r(n) b(n-1).
\end{align*}
\end{proof}
\end{lem}

The above estimate is sharp.  For any $n > 1$, let $r(1) = \cdots
= r(n-1)  = 1$ and $r(n) = n$.  Then $b(1) = \cdots = b(n-1) = 1$,
and $b(n) = n $, so $b(n)/b(n-1) = n = r(n)$.

\begin{cor} \label{cor: b(n)+1/b(n) <= prod}
For all $n, m \geq 1 $
\[
    \frac{b(n+m)}{b(n)} \leq \prod_{k= n+1}^{n+m} r(k).
\]
\begin{proof}
Write ${b({n+m})}/{b(n)}$ as a telescoping product, and apply
Lemma \ref{lem: b(n) <= r(n)} $m$ times:
\[
    \frac{b({n+m})}{b(n)} =  \prod_{k= n+1}^{n+m}  \frac{b(k)}{b({k-1})}
    \leq \prod_{k= n+1}^{n+m} r(k).
\]
\end{proof}
\end{cor}

From the above estimate, it follows that $b(n) \leq
\prod_{k=1}^{n} r(k)$. Which implies $b(n) \leq n!$.  However,
from Lemma \ref{lem: max growth 2 n-1} we know in fact that $b(n)
\leq 2^{n-1}$.  So while Lemma \ref{lem: b(n) <= r(n)} gives a
sharp estimate of the short term growth of $b(n)$, in the long
term it is highly inaccurate. However, we only use the above
corollary to obtain lower bounds on growth, so the inaccuracy is
somewhat reduced. We state the reciprocal of it for reference.

\begin{cor} \label{cor: b(n)/b(n)+m => prod}
For all $n, m \geq 1 $
\[
    \frac{b(n)}{b({n+m})} \geq \prod_{k= n+1}^{n+m} 1/r(k).
\]
\end{cor}

We need an estimate on the ratio of sums.  The proof is trivial.

\begin{lem} \label{lem: sum a_i /  sum c_k}
Let $\alpha_1, \dots, \alpha_l$ and $\beta_1, \dots, \beta_m$ be
positive real numbers. If for all $i$ and $j$ we have $\alpha_i
\leq \beta_j$, then
\[
    \frac{\sum_{i=1}^{l} \alpha_i}{\sum_{j=1}^m \beta_j} \leq \frac{l}{m}.
\]
\end{lem}

\begin{lem} \label{lem: b(n)-r/b(n)-s for r> s}
If $r > s > 1$, then for any $n \geq r$
\[
    \frac{\sum_{k=1}^{r} b(n-k)}
    {\sum_{k=1}^{s} b(n-k)} \leq \frac{r}{s}.
\]
\begin{proof}
We have
\[
     \frac{\sum_{k=1}^{r} b(n-k)}
 {\sum_{k=1}^{s} b(n-k) }=
     \frac{\sum_{k=1}^{s} b(n-k) }{\sum_{k=1}^{s} b(n-k) }
    + \frac{\sum_{k=s+1}^{r} b(n-k) }
    {\sum_{k=1}^{s} b(n-k) }.\\
\]
{Since the $b(n)$ are non-decreasing, Lemma \ref{lem: sum a_i /
sum c_k} applies to the second term, and}
\[
    \frac{\sum_{k=1}^{r} b(n-k) }
   {\sum_{k=1}^{s} b(n-k) }
    \leq 1 + \frac{r-s}{s} = \frac{r}{s}.
\]

\end{proof}
\end{lem}

When there are more terms in the denominator, we obtain the
following corollary in a similar fashion.
\begin{cor} \label{cor: b(n)-r/b(n)-s for r< s}
If $s > r > 1$, then for any $n \geq s$
\[
    \frac{\sum_{k=1}^{r} b(n-k) }
    {\sum_{k=1}^{s} b(n-k) }\geq \frac{s}{r}.
\]
\end{cor}

Recall that $\lambda(n) = 1 + [r(n) - 1]/r(n-1)$.

\begin{lem} \label{lem: b(n)+1/b(n) <= 1 + lambda(n+1)}
For any $n \geq 0$, if $\Delta r(n+1) > 1 $, then
\[
    \frac{b(n+1)}{b(n)} \leq \lambda(n+1).
\]
\begin{proof}
\begin{align*}
    \frac{b(n+1)}{b(n)} &=
    \frac{b(n) + b(n-1) + \dots + b({n+1 -r(n+1)})}{b(n)}\\
    &= 1 + \frac{b(n-1) + \dots + b({n+1 -r(n+1)})}{b(n-1) + \dots +
    b(n-r(n))}.
\end{align*}
Note that there are $r(n+1) -1$ terms in the numerator, $r(n)$
terms in the denominator, and by assumption $r(n+1) -1
> r(n)$. Thus, we can use Lemma \ref{lem:
b(n)-r/b(n)-s for r> s}  to obtain
\[
     \frac{b(n+1)}{b(n)} \leq 1 + \frac{r(n+1) - 1}{r(n)} =
     \lambda(n+1).
\]
\end{proof}
\end{lem}

\begin{lem} \label{lem: b(n)+1/b(n) => 1 + lambda(n+1)}
For any $n \geq 0$, if $\Delta r(n+1) <  1$, then
\[
    \frac{b(n+1)}{b(n)} \geq \lambda(n+1).
\]
\begin{proof}
Similar to Lemma \ref{lem: b(n)+1/b(n) <= 1 + lambda(n+1)}, except
that we use Corollary \ref{cor: b(n)-r/b(n)-s for r< s}.
\end{proof}
\end{lem}

\begin{rem}
For the $r$-generalized Fibonacci numbers $(f_{r,n})$, this
estimate shows that $f_{r, n+1}/f_{r,n}  \geq 1 + (r-1)/r = 2
-1/r$ for $n
> 2r - 1$.
\end{rem}

Notice that $\lambda(n)$ is either an upper bound or a lower bound
depending on $\Delta r(n)$.

The above lemma gives us information about the asymptotics of
$b(n)$.

\begin{cor} \label{cor: Omega(1+(m-1)/M}
If $m = \liminf_{n \to \infty} r(n)$ and $M = \limsup_{n \to
\infty} r(n)$, then $b(n)$ is $\Omega(1+\frac{m-1}{M})$.
\end{cor}

Hence, if $\liminf r(k) > 1$ and $\limsup r(k) < \infty$ , then
the $b(n)$ grow exponentially fast. In contrast, the growth rate
of many meta-Fibonacci sequences is of only linear order. For
instance, the Conway sequence
\[
    a(n) = a(a(n-1)) + a(n - a(n-1)), \quad n \geq 3,
\]
$a(1) = a(2) =1$.  It is known that $ \lim_{n \to \infty} a(n)/n =
1/2$ \cite{Mal91}.  We discuss this phenomenon in more detail in
$\S$\ref{sect: Asymptotics}.

Recall that $\mu(n,s) = 2 + [\Delta r(n) - 1] \prod_{k= n -
s}^{n-1} 1/r(k)$.  We give estimates on growth, in terms of
$\mu(n,s)$.

\begin{lem} \label{lem: b(n)/b(n)-1 <= mu(n, r(n-1))}
For any $n \geq 1$, if $\Delta r(n)  < 1$, then
\[
    \frac{b(n)}{b(n-1)} \leq  \mu(n, r(n-1)).
\]
\begin{proof}
Using Lemma \ref{lem: Delta r(n) = 1 -> growth =2} we have

\begin{align*}
    &\frac{\sum_{k=1}^{r(n-1)+1} b(n-k) }{b(n-1)} =2\\
     &\frac{\sum_{k=1}^{r(n)} b(n-k) }{b(n-1)} +
     \frac{\sum_{k=r(n) + 1}^{r(n-1)+1} b(n-k) }{b(n-1)}
     = 2\\
    &\frac{b(n)}{b(n-1)} =
    2 - \frac{\sum_{k=  r(n) + 1}^{r(n-1)+1} b(n-k) }{b(n-1)}\\
    &\frac{b(n)}{b(n-1)} \leq 2 - [r(n-1) -r(n) + 1]
    \frac{b({n-r(n-1)-1})}{b(n-1)}\\
    \intertext{since the $b(n)$ are non-decreasing, so}
    &\frac{b(n)}{b(n-1)} \leq 2 - \{-[r(n) - r(n-1) - 1] \} \prod_{k=n-r(n-1)}^{n-1}1/r(k)\\
\end{align*}
  by Corollary \ref{cor: b(n)/b(n)+m => prod}. The right-hand side
is $\mu(n,r(n- 1))$, so the lemma is shown.
\end{proof}
\end{lem}

\begin{rem}
For the $r$-generalized Fibonacci numbers $(f_{r,n})$, this
proposition implies that $ f_{r,n+1}/f_{r,n} \leq 2 - r^{-r}$, for
$n > 2r - 1$.
\end{rem}

\begin{lem} \label{lem: b(n)/b(n)-1 => mu(n, r(n-1)+1)}
If $\Delta r(n)  > 1$, then
\[
    \frac{b(n)}{b(n-1)} \geq  \mu(n, r(n)-1).
\]
\begin{proof}
By Lemma \ref{lem: Delta r(n) = 1 -> growth =2} we have
\[
    \frac{\sum_{k=1}^{r(n-1)+1} b(n-k) }{b(n-1)} =2.
\]
Thus,
\begin{align*}
    \frac{b(n)}{b(n-1)} &= \frac{\sum_{k=1}^{r(n)} b(n-k)}{b(n-1)}\\
           &= \frac{\sum_{k=1}^{r(n-1)+1} b(n-k) }{b(n-1)} +
         \frac{\sum_{k=r(n-1)+2}^{r(n)} b(n-k) }{b(n-1)}\\
         &= 2+\frac{\sum_{k=r(n-1)+2}^{r(n)} b(n-k) }{b(n-1)}\\
            &\geq 2 + [r(n) -r(n-1) - 1] \frac{b(n-r(n))}{b(n-1)}\\
         \intertext{since the $b(n)$ are non-decreasing, }
          &\geq 2 + [r(n) - r(n-1) - 1]\prod_{k=n-r(n)+1}^{n-1} 1/r(k)\\
\end{align*}
by Corollary \ref{cor: b(n)/b(n)+m => prod}. The right-hand side
is $\mu(n,r(n)- 1)$, so the lemma is shown.

\end{proof}
\end{lem}

As with $\lambda(n)$, we have $\mu(n,R)$, where $R =  \max
\set{r(n) - 1, r(n-1)}$, is either an upper and lower bound for
growth depending on $\Delta r(n)$.

We now compare $\mu$ to $\lambda$.

\begin{lem} \label{lem: mu => lambda}
If $\Delta r(n)  \leq 1$ for some $n \geq 1$, then
\[
    \mu(n, r(n) - 1) \geq \lambda(n).
\]
\begin{proof}
We have $[r(n) - r(n-1) - 1 ] \leq 0 $ by assumption.  Also
$\prod_{k = n-r(n) + 1}^{n-2} 1/r(k) \leq 1$, since $r(k) \geq 1$
for all $k$.  Thus,
\[
    [r(n) - r(n-1) - 1 ]
    \left[  -1 + \prod_{k = n-r(n) + 1}^{n-2}  1/r(k)  \right]
    \geq 0.
\]
The lemma follows by straightforward algebra.
\end{proof}

\end{lem}

\begin{lem} \label{lem: mu <= lambda}
If $\Delta r(n) \geq 1$ for some $n \geq 1$, then
\[
    \mu(n, r(n - 1)) \leq \lambda(n).
\]
\begin{proof}
Similar to the above lemma.
\end{proof}

\end{lem}

We are now ready to prove the Main Theorem.  We consider various
cases for $\Delta r(n)$.  We then combine appropriate estimates of
$b(n)/b(n-1)$ in terms of $\mu$ and $\lambda$.

\begin{proof}[Proof of Main Theorem.]
Fix $n \geq 1$.  If $\Delta r(n)  = 1$, then by Lemma \ref{lem:
Delta r(n) = 1 -> growth =2} we know $b(n)/b(n-1) =2 = \lambda(n)
= \mu(n, \cdot)$ and we are done.

If $\Delta r(n)  < 1$, then
\[
    \lambda(n) \leq \frac{b(n)}{b(n-1)} \leq \mu(n, r(n-1)),
\]
with the first inequality by Lemma \ref{lem: b(n)+1/b(n) => 1 +
lambda(n+1)} and the second by Lemma \ref{lem: b(n)/b(n)-1 <=
mu(n, r(n-1))}. Additionally, $\lambda(n) \leq  \mu(n, r(n) - 1) $
by Lemma \ref{lem: mu => lambda}. Hence
\[
    \min \set{\lambda(n), \mu(n, r(n) - 1)} \leq \frac{b(n)}{b(n-1)}.
\]

 Finally, if $\Delta r(n)  > 1$, then
\[
   \mu(n, r(n)-1) \leq \frac{b(n)}{b(n-1)} \leq  \lambda(n) ,
\]
by Lemma \ref{lem: b(n)/b(n)-1 => mu(n, r(n-1)+1)} and Lemma
\ref{lem: b(n)+1/b(n) <= 1 + lambda(n+1)}.  Also, $\mu(n, r(n-1))
\leq \lambda(n) $ by Lemma \ref{lem: mu => lambda}, so
\[
   \frac{b(n)}{b(n-1)} \leq \max \set{\lambda(n),
    \mu(n, r(n-1))}.
\]

Therefore, we can combine the above inequalities in all cases to
obtain:
\[
    \min \set{\lambda(n), \mu(n, r(n)-1)} \leq  \frac{b(n)}{b(n-1)} \leq \max \set{\lambda(n),
    \mu(n, r(n-1))}.
\]
\end{proof}

\section{Asymptotic Growth} \label{sect: Asymptotics}

In this section we examine the asymptotic growth of $b(n)$. We
compare the growth rate of $b(n)$ to other families of
meta-Fibonacci sequences, which is polynomial order in all known
cases. We show that $b(n)$ can have a variety of different growth
rates: exponential, linear, and logarithmic. However, the possible
asymptotic limits for $b(n)$ are restricted.

To date two other families of meta-Fibonacci sequences have
appeared in print, see \cite{DGNW} and \cite{CCT}.  Sub-families
of the latter family are also studied in \cite{HT93b} and
\cite{JR05}.

In \cite{DGNW} $(p,q)$-sequences were introduced.  A $(p,q)$
sequence $(F_n)$ is defined as follows.  For fixed positive
integers $p$ and $q$, and values $a_1,\cdots,a_p$, let $F_n=a_n$
with probability one for $n\leq p$ and set
$F_{n+1}=\sum^q_{k=1}F_{j_k}$ for $n\geq p$, where the $j_k$ are
randomly chosen, with replacement, from $(1,2,\cdots,n)$.  They
give only probabilistic results. They show that the expected value
of $F_n$ grows as a polynomial in $n$ of degree $q-1$ \cite[Thm.
1]{DGNW}.

In \cite{CCT} J. Callaghan, J. Chew and S. Tanny studied a family
of sequences parameterized by $a > 0$, $ k >1$:
\[
    T_{a,k}(n) = \sum_{i=0}^{k-1} T_{a,k}(n  - i - a
    -T_{a,k}(n-i-1)), \quad n > a+k, \ k \geq 2
\]
with $ T_{a,k}(n) = 1 $ for $1 \leq n \leq a+k$.  For $k$ odd, the
growth rate of sequences in this family is linear \cite[Cor.
5.14]{CCT}; for all $a$ and all odd $k$
\[
    \lim_{n \to \infty } \frac{ T_{a,k}(n)}{n} = \frac{k-1}{k}.
\]

We now examine the asymptotics of variable-$r$-bonacci sequences.
In contrast to the above meta-Fibonacci sequences, but like the
$r$-generalized Fibonacci sequences, $b(n)$ can grow
exponentially. As previously noted, by Lemma \ref{lem: max growth
2 n-1}, all such sequences are $O(2^{n-1})$. By Corollary
\ref{cor: Omega(1+(m-1)/M}, t$\Omega([1+(m-1)/M]^n)$ for any $ 2
\leq m \leq M $.  It is possible that that $b(n) \sim \gamma^n$,
but the possible values of $\gamma$ are limited

\begin{lem} \label{lem: lim inf growth <= 2}
For any variable-$r$-bonacci sequence $b(n)$, we have
\[
    \liminf_{n \to \infty} \frac{b(n)}{b(n-1)} \leq 2.
\]
 \begin{proof}
Suppose not. By Theorem \ref{thm: b(n) growth}.d, for all $n$
sufficiently large $\Delta r(n) > 1$.  It follows that for $n$
large $n -r(n) < (n-1) -r(n-1)$. Thus, for $n$ large $(n -r(n))$
is a strictly decreasing sequence of integers. Therefore, $N- r(N)
< 0$ for some $N$.  Contrary to  $r(n) \leq n$ by Definition
\ref{defn: vrmFs}.
 \end{proof}
\end{lem}

The possible limits for the sequence $(b(n)/ b(n-1))$ are
restricted.

\begin{lem}
If $\lim_{n \to \infty} r(n) = R$, then
\[
    \lim_{n \to \infty} \frac{b(n)}{ b(n-1)} = \alpha_R,
\]
where $\alpha_1 =1$.
\begin{proof}
For $n$ sufficiently large, $b(n)$ satisfies the $R$-boncacci
recursion or is eventually constant.
\end{proof}
\end{lem}

The only other possible limit of $b(n)/b(n-1)$ is 2.

\begin{prop}
If the sequence $b(n)/ b(n-1)$ converges and $r(n)$ is not
eventually constant, then
\[
    \lim_{n \to \infty} b(n)/{ b(n-1)} = 2.
\]
\begin{proof}
If $ \limsup_{n \to \infty} {b(n)}/{b(n-1)} < 2$, then by
Corollary \ref{cor: lim sup growth < 2 -> r(n) const}, $\lim_{n
\to \infty} r(n) = R$ for some $ R \in \bZ^+$, contrary to
assumption. Thus $ \limsup_{n \to \infty} {b(n)}/{b(n-1)} \geq 2$.
By Lemma \ref{lem: lim inf growth <= 2}, $ \liminf_{n \to \infty}
{b(n)}/{b(n-1)} \leq 2.$ Therefore, the only possible limit is 2.

\end{proof}
\end{prop}

\begin{cor}\label{cor: asymptotic grwoth}
If $\lim_{n\to \infty}b(n)/\gamma^n $ exists for some $\gamma \in
\bR$, then $\gamma = \alpha_R$ for some $R \geq 1$, or $\gamma
=2$.
\end{cor}

From Example \ref{ex: b(n) max growth}, we know that $  \lim_{n
\to \infty} b(n)/{ b(n-1)} = 2$ occurs.  The following example
shows that it occurs in a non-trivial case.

\begin{example}
For $n \geq 2$, let $r(n) = n$ for $n$ even, and $r(n) = n-1$ for
$n$ odd.
\[
\begin{array}{|c|c|c|c|c|c|c|c|c|c|c|}\hline
  n    & 0 & 1  &  2 &  3 & 4  & 5  &  6 & 7  & 8  & 9  \\
  \hline
  r(n) & 1 & 1  &  2 &  2 & 4  & 4  &  6 & 6  & 8  & 8  \\
  \hline
  b(n)  & 1 & 1  &  2 &  3 &  7 & 13  & 27  & 53  & 107  & 213  \\
  \hline
\end{array}
\]
We claim that if $n > 2$, then $b(n) = 2 b(n-1) +1$ for $n$ even,
and $b(n) = 2 b(n-1) - 1$ for $n$ odd.  Hence, $ \lim_{n \to
\infty} b(n)/{ b(n-1)} = 2$.  The proof is left as an exercise.
\end{example}

The following example shows that the sequence $b(n)/ b(n-1)$ need
not converge.

\begin{example}
For $n \geq 2$, let $r(n) = 2$ for $n$ even, and $r(n) = 3$ for
$n$ odd.
\[
\begin{array}{|c|c|c|c|c|c|c|c|c|c|c|}\hline
  n    & 0 & 1  &  2 &  3 & 4  & 5  &  6 & 7  & 8  & 9  \\
  \hline
  r(n) & 1 & 1  &  2 &  3 & 2  & 3  &  2 & 3  & 2  & 3  \\
  \hline
  b(n)  & 1 & 1  &  2 &  4 &  6 & 12  & 18  & 36  & 54  & 108  \\
  \hline
\end{array}
\]
It is left as an exercise to show that $b(n)/ b(n-1) = 2$ for $n
> 2$ and odd, and $b(n)/ b(n-1) = 3/2$ for $n>2$ and even.
\end{example}

By taking $r(n) =1 $ fairly often, we can have linear growth for
$b(n)$.

\begin{example}
For $n \geq 2$, let $r(n) = 2$ if $n = 2^k$ for some $k \in \bZ$,
and $r(n) = 1$ otherwise.
\[
\begin{array}{|c|c|c|c|c|c|c|c|c|c|c|}\hline
  n    & 0 & 1  &  2 &  3 & 4  & 5  &  6 & 7  & 8  & 9  \\
  \hline
  r(n) & 1 & 1  &  2 &  1 & 2  & 1  &  1 & 1  & 2  & 1  \\
  \hline
  b(n)  & 1 & 1  &  2 &  2 &  4 & 4  & 4  & 4  & 8  & 8  \\
  \hline
\end{array}
\]
It is easy to show that $n/2 \leq b(n) \leq n$ for $n \geq 1$.
That is, $b(n)$ is $\Theta(n)$.
\end{example}

However in the case of linear growth, we cannot have an asymptotic
limit other than zero.

\begin{prop}\label{prop: lin growth limit}
    If $\lim_{n \to \infty} b(n)/n = L$, where $0 \leq L <
    \infty$, then $L = 0$.
\begin{proof}
Contrarily,  if $L > 0$ we can take  $0 < \varepsilon \ll L$. We
can then find some large $N$ such that
\begin{enumerate}
    \item $r(N) = 1$ (or else the growth is exponential);
    \item $r(N+1) > 1$ (or else $b(n)$ is eventually constant);
    \item $N/(N+1) > 1 - \varepsilon$;
    \item for all $n \geq N$, $|b(n)/n - L|< \varepsilon$.
\end{enumerate}
We have
\begin{align*}
      \frac{b(N+1)}{N+1}   &= \frac{b(N+1)}{b(N)} \frac{b(N)}{N}
      \frac{N}{N+1}\\
      &\geq(2)(L-\varepsilon)(1 - \varepsilon),\\
      &= 2L + O(\varepsilon).
\end{align*}
by Lemma \ref{lem: Delta r(n) = 1 -> growth =2}, condition 4, and
condition 3 respectively.  So,
\[
    \left| \frac{b(N+1)}{N+1} - L \right| \geq 2L + O(\varepsilon) - L
    = L + O(\varepsilon)  .
\]
But by condition 4, we have
\[
    \left| \frac{b(N+1)}{N+1} - L \right| < \varepsilon \ll L.
\]
Therefore, $L = 0$.
\end{proof}

\end{prop}

We can have slower than polynomial growth.  No other known
Fibonacci-type sequence grows so slowly.

\begin{example}\label{eg: log growth}
For $n \geq 2$, let $r(n) = 2$ if $n = 2^{2^k}$ for some $k \in
\bZ$, and $r(n) = 1$ otherwise.
\[
\begin{array}{|c|c|c|c|c|c|c|}\hline
  n    & 0 & 1  &  2 &  4 & 16  & 256  \\
  \hline
  r(n) & 1 & 1  &  2 &  2 & 2  & 2  \\
  \hline
  b(n)  & 1 & 1  &  2 &  4 &  8 & 16    \\
  \hline
\end{array}
\]
It is easy to show that $b(n)$ is $\Theta(\log_2 n)$.
\end{example}

Similarly, we can construct examples that are $\Theta(\log_2
\log_2 n)$, etc.  Thus $b(n)$ can grow quite slowly indeed.

\section{Generalization} \label{sect: Generalization}
We define a generalization of $b(n)$.  This generalization allows
us to pick different initial conditions for our sequence.  It also
allows us to remove the restrictions that $r$ be sublinear.

\begin{defn}\label{defn: ex vrmFs}
We call a double sequence $\beta(n), \ n \in \bZ$, an
\emph{extended variable-$r$ meta-Fibonacci sequence} if there
exists $r: \bZ \to \bZ^+$ such that for all $n \in \bZ$
\[
    \beta(n) = \sum_{k=1}^{r(n)} \beta(n-k).
\]
\end{defn}
Note that there no restriction that $r$ be sublinear.  Provided
$\beta(n) > 0$ for all $n \in \bZ$, all results in this paper
apply to an extended $\beta(n)$, except Lemma \ref{lem: max growth
2 n-1}. Similarly, if $\beta(n) < 0$  for all $n \in \bZ$, all
results in this paper, except Lemma \ref{lem: max growth 2 n-1},
are easily generalized.  The behavior of $\beta(n)$ with both
positive and negative terms is an interesting question.

Given $r: \bN \to \bZ^+$, we can define a sequence $\beta(n)$
generated by $r$ as follows.  Pick any real number $\beta(-1)$ as
an initial condition. For $n \leq -1$, let $\beta(n) = \beta(-1)$
and let $r(n) = 1$. For $n \geq 0$, define $\beta(n)$ by the
variable-$r$-bonacci recursion.

Now let $r: \bN \to \bZ^+$ and let $M_r = \sup_{n \in \bN} r(n) -n
$. We consider $r$ with $M_r$ finite.  Note that $r(0) - 0 > 0$,
so $M_r \geq 1$.   We give a construction for extending $r$ to a
function on all integers, so that it generates an extended
variable-$r$-bonacci sequence.

\begin{defn}\label{defn: ex vrmFs recipie}
Let $r: \bN \to \bZ^+$ and let $M_r = \sup_{n \in \bN} r(n) -n <
\infty$.  Choose $\beta(-1), \dots, \beta(-M_r) \in \bR$.  For $n
\geq 0$ let
\[
    \beta(n) = \sum_{k=1}^{r(n)} \beta(n-k).
\]
 For $n =  -1, -2, \dots$ let $r(n) = M_r$, and let
\[
    \beta(n - M_r) = \beta(n) - \sum_{k=1}^{M_r - 1} \beta(n-k).
\]
\end{defn}

\begin{prop}
Let $r: \bN \to \bZ^+$ with  $M_r = \sup_{n \in \bN} r(n) -n <
\infty$. The double sequence $\beta(n)$ constructed in Definition
\ref{defn: ex vrmFs recipie
} is an extended variable-$r$
meta-Fibonacci sequence generated by the extension of $r(n)$ to
$\bZ$.
\begin{proof}
It is straightforward to check that $\beta(n)$ is well defined for
all $n \in \bZ$,
 and it satisfies the correct recursion
relation.
\end{proof}
\end{prop}

\begin{cor}
Let $r: \bN \to \bZ^+$ with  $M_r = \sup_{n \in \bN} r(n) -n <
\infty$.  The set of variable-$r$ meta-Fibonacci sequences
generated by $r(n)$ is an $M_r$-dimensional real vector space.
\end{cor}

In particular, if $r$ is sublinear, $M_r = 1$, so we can use this
construction to get a sequence $\beta(n)$ depending on one initial
condition $\beta(-1)$.

\providecommand{\bysame}{\leavevmode\hbox
to3em{\hrulefill}\thinspace}

\providecommand{\href}[2]{#2}

AMS Subject Classification: {11B37, 11B39, 11B99.}


\begin{thebibliography}{DGNW}

\bibitem[BH]{BH92}
Bodil Branner and John Hubbard, \emph{Iteration of cubic
polynomials, part
  {I}{I}: Patterns and parapatterns}, Acta Math. \textbf{169} (1992), 229--
  325.

\bibitem[CCT]{CCT}
Joseph Callaghan, John~J. Chew, III, and Stephen~M. Tanny,
\emph{On the
  behavior of a family of meta-{F}ibonacci sequences}, SIAM J. Discrete Math.
  \textbf{18} (2004), no.~4, 794--824.

\bibitem[Co]{Cnly89}
{B}.~{W}. Conolly, \emph{Meta-{F}ibonacci sequences}, ch.~XII of
\emph{Fibonacci \& {L}ucas {N}umbers and the {G}olden {S}ection},
by S. Vajda,
 pp.~127--138, {E}llis {H}orwood, 1989.

\bibitem[DGNW]{DGNW}
R.~Dawson, G.~Gabor, R.~Nowakowski, and D.~Wiens, \emph{Random
{F}ibonacci-type
  sequences}, Fibonacci Quart. \textbf{23} (1985), no.~2, 169--176.

\bibitem[E1]{E01}
Nathaniel~D. Emerson, \emph{Dynamics of polynomials whose {J}ulia
set is an area zero
  {C}antor set}, Ph. {D}. thesis, University of California Los Angeles, 2001.

\bibitem[E2]{E05}
\bysame, \emph{Return times of polynomials as meta-{F}ibonacci
  numbers}, Preprint, August 2005.



\bibitem[Ho]{Hof79}
Douglas~R. Hofstader, \emph{G{\"o}del, {E}scher, {B}ach: an
{E}ternal {G}olden
  {B}raid}, Basic Books, 1979.

\bibitem[HT]{HT93b}
J.~Higham and S.~Tanny, \emph{More well-behaved meta-{F}ibonacci
sequences},
  Proceedings of the Twenty-fourth Southeastern International Conference on
  Combinatorics, Graph Theory, and Computing (Boca Raton, FL, 1993), vol.~98,
  1993, pp.~3--17.

\bibitem[JR]{JR05}
Brad Jackson and Frank Ruskey, \emph{Meta-{F}ibonacci sequences,
binary trees,
  and extremal compact codes}, Preprint, April 2005.

\bibitem[Ma]{Mal91}
Colin~L. Mallows, \emph{Conway's challenge sequence}, Amer. Math.
Monthly
  \textbf{98} (1991), no.~1, 5--20.

\bibitem[Mi]{Miles}
E.~P. {Miles, Jr.}, \emph{Generalized {F}ibonacci numbers and
associated
  matrices}, Amer. Math. Monthly \textbf{67} (1960), no.~8, 745--752.

\end{thebibliography}
\end{document}